\documentclass[11pt ]{article}
\newtheorem{definition}{Definition}
\newtheorem{lemma}[definition]{Lemma}
\newtheorem{proposition}[definition]{Proposition}
\newtheorem{corollary}[definition]{Corollary}
\newtheorem{theorem}[definition]{Theorem}

\begin{document}
\centerline{\bf The Bernstein-Gelfand-Gelfand complex for the groups $Sp(n,1)$ and $F_{4(-20)}$}
\centerline {\bf as a Kasparov  module}

\bigskip\bigskip
\centerline{\bf Pierre Julg}
\bigskip
\centerline {Universit\'e d'Orl\'eans}
\bigskip \bigskip\bigskip\bigskip\bigskip

\noindent{\bf 0. Introduction.}

\bigskip
Let $G$ be a locally compact group and $A$ be any $G-C^*$-algebra. We say that $G$ satisfies the Baum-Connes
conjecture with coefficients in $A$ if that the Baum-Connes assembly map:
$$K_*^G({\bf E}G,A)\rightarrow K_*(C^*_r(G,A))$$
is an isomorphism.
The left hand side is the $G$-equivariant $K$-homology of the universal proper space for
proper $G$-actions, which is well defined up to $G$-equivariant homotopy, with coefficients in
$A$. The right hand side is the
$K$-theory of the reduced crossed product of $A$ by $G$.

We say that $G$ satisfies the Baum-Connes conjecture with coefficients if it satisfies the conjecture with
coefficients in any $G-C^*$-algebra $A$. On the other hand, $G$ satisfies the usual Baum-Connes conjecture, i.e. the
conjecture without coefficents, if it satisfies the conjecture with scalar coefficients, i.e. for 
$A={\bf C}$, the algebra of complex numbers. The important point about the conjecture {\it  with} coefficients is
that it is stable by restriction to closed subgroup. In particular, if $G$ is a group which satisfies the
Baum-Connes conjecture with coefficients, then any closed, and in particular any discrete subgroup of $G$ satisfies
the Baum-Connes conjecture.

We are interested in the following situation: let us assume that $G$ is a semi-simple Lie group, connected
with finite centre. The conjecture without coefficients for $G$ is known to be true. There are actually two
completely distinct proofs of that fact. In 1984, A. Wassermann [W](following the work of Pennington-Plymen and
Valette) has proved the conjecture using the Harish-Chandra theory which allows to describe the reduced dual
of such groups. However the philosophy of A. Connes non commutative geometry suggested that one should provide a
proof of the conjecture completely independent of such a representation theory. This has been done in the
remarkable work of V. Lafforgue in 1998 [L].  

On the other hand, let $\Gamma$ be a discrete subgroup of $G$, the question of Baum-Connes for $\Gamma$ can be
summerized as follows: 

1) injectivity is known to be true since the work of Kasparov on the Novikov conjecture in
1980 [K1][K3][; in fact this is a consequence of the injectivity of the Baum-connes assembly map with coefficient in
any
$A$ for the Lie group $G$.

2) the question of surjectivity of the Baum-Connes assembly map for the discrete group $\Gamma$, or more generally
the surjectivity of the Baum-Connes assembly map with coefficients in any $A$ for the Lie group $G$ is a difficult
problem and a test for the conjecture. Let us list a few results and open problems. 

In the case where $G$ is simple with real rank one, there is a dichotomy between two classes of groups:

(i) if $G$ is (locally isomorphic to) one of the groups $S0_0(n,1)$ or $SU(n,1)$ ($n\geq 2$), $G$ has the Haagerup
property. Therefore, by the Higson-Kasparov theorem [HK][J4], $G$ satisfies the Baum-Connes conjecture with
coefficients, and so do all its discrete subgroups. Moreover, $G$ is $K$-amenable [JV]: the reduced crossed products
$C^*_r(G,A)$ can be replaced  by the full crossed products $C^*(G,A)$. 

The case of
$SO_0(n,1)$ had  in fact been  solved by Kasparov [K2] already in 1983, combining geometry (the conformal structure
on the sphere at infinity of the hyperbolic space) and the theory of unitary representations (complementary series).
This  was generalized to
$SU(n,1)$ by Kasparov and the author [JK].

(ii) if $G$ is (locally isomorphic to) one of the groups $Sp(n,1)$ ($n\geq 2$) or $F_{4(-20)}$ , then $G$ has
Kazhdan's property $T$. This fact makes the Baum-Connes conjecture more difficult since the full and reduced
crossed product do not have in general the same $K$-theory. The first deep result in that direction was obtained by
V. Lafforgue in 1998 [L]: if $\Gamma$ is a cocompact discrete subgroup of such a group $G$, then $\Gamma$ satisfies
the Baum-connes conjecture (without coefficients), and even the conjecture with
coefficients in a commutative $C^*$-algebra.

In this paper we present the main step in the proof of  the following result:

\begin{theorem}
 Let $G$ be one of the groups $Sp(n,1)$ ($n\geq 2$) or $F_{4(-20)}$ . Then $G$ satisfies the
Baum-Connes conjecture with coefficients. In particular all discrete subgroups of $G$ satisfy the Baum-Connes
conjecture.
\end{theorem}

Our proof is a generalisation of the method of [JK]. The difference is that it is no longer possible to use  the
theory of unitary representations since the complementary series stays away from the trivial representation. The
main ingredient in our proof is the use of a family of uniformly bounded representations , constructed by M.
Cowling, which approach the trivial representation.

Recently, V. Lafforgue has shown that any Gromov hyperbolic group $\Gamma$ satisfies the Baum-Connes conjecture with
coefficients. Note that the range of application of Lafforgue's method and ours are in general position: the case
where
$\Gamma$ is a cocompact discrete subgroup of
$Sp(n,1)$ ($n\geq 2$) or
$F_{4(-20)}$ can therefore be obtained by both. However we obtain the case of all discrete subgroups of such Lie
groups, and Lafforgue obtains all Gromov hyperbolic groups.

Let us say a few words about the case of higher rank. Very few is known. Let $G$ be a simple Lie group of real rank
at least 2, and $\Gamma$ a discrete subgroup of $G$.

(i) Assume that  $G$ is  
$SL_3({\bf R})$ or
$SL_3({\bf C})$.
 V. Lafforgue has shown  that any discrete cocompact subgroup
$\Gamma$ of $G$ satisfies the Baum-connes conjecture (without coefficients). I. Chatterji has generalized this fact
to $SL_3({\bf H})$ and $E_{6(-26)}$. It is not known whether such a group
$\Gamma$ (or {\it a fortiori} the  Lie group $G$) satisfies the conjecture with coefficients. 
Moreover, nothing is known about the Baum-Connes conjecture for general discrete
subgroups of $G$.

(ii) If $G$ is either another real rank 2 simple group, or a simple group with real rank at least 3, nothing is
known.

In particular there are two very difficult open problems: Prove that $SL_3({\bf Z})$ satisfies the Baum-Connes conjecture, or prove that
$SL_3({\bf R})$ satisfies the Baum-Connes with coefficients.

On the other hand, leaving the world of Lie groups, the so called Gromov random groups (or Gromov's monsters) are
known to be counterexamples to the Baum-Connes conjecture with coefficients (but possibly not to the plain
conjecture).

\bigskip\bigskip

\noindent{\bf 1.$K$-theory and Fredholm modules.}

\bigskip
\noindent{\bf 1.1 Review of Kasparov's  Dirac-dual Dirac method.}
\bigskip
\noindent

 From now on,
$G$ will be a semi-simple Lie group, connected with finite centre. Let
$K$ be a maximal compact subgroup of $G$. The Riemannian symmetric space
$X=G/K$ representative of the universal space
${\bf E}G)$ . The subgroup $K$ is the stabilizer of a point $x_0$ in the symmetric space and
therefore acts by a linear representation in the real vector space $V$ tangent to $X$ at
point $x_0$. The left hand side is isomorphic to  the
$K$-theory of the crossed product of $C_0(V^*)\otimes A$ by $K$. The Baum-Connes conjecture for that case is
 the so-called Connes-Kasparov conjecture: the map
$$\alpha_A: K_*(C^*(K,C_0(V^*)\otimes A))\rightarrow K_*(C^*_r(G,A))$$ is an isomorphism.
That map $\alpha_A$ is called the Dirac map because it can be constructed using a Dirac
operator on $G/K$.

G. Kasparov (cf [K1][K3]) has shown that the map $\alpha_A$ is injective. More precisely he constructed a left
inverse $\beta_A : K_*(C^*_r(G,A))\rightarrow K_*(C^*(K,C_0(V^*)\otimes A))$, called dual
Dirac. The composition $\beta_A\alpha_A$ is the identity of $K_*(C^*(K,C_0(V^*)\otimes A))$,
but the composition $\gamma_A=\alpha_A\beta_A$ is {\it a priori} only an idempotent map of the
$K$-theory group $K_*(C^*_r(G,A))$.
The conjecture is therefore equivalent to the fact  that the map $\gamma_A$ is the identity.

The maps $\alpha_A$ and $\beta_A$ are described by Kasparov thanks to explicit elements in equivariant 
$KK$-theory:
$$\alpha\in KK_G(C_0(T^*X),{\bf C})$$
$$\beta\in KK_G({\bf C}, C_0(T^*X))$$
such that $\alpha\otimes\beta =1$ in $KK_G(C_0(T^*X),C_0(T^*X))$.
 In particular , the map $\gamma_A$ comes from the
 idempotent element
$\gamma=\beta\otimes\alpha$ of the commutative ring $R(G)=KK_G({\bf C},{\bf C})$. The
abelian groups
$K_*(C^*_r(G,A))$ are modules over the ring $R(G)$ via the following ring homomorphisms:
$$R(G)=KK_G({\bf C},{\bf C})\rightarrow KK_G(A,A)\rightarrow KK(C^*_r(G,A),C^*_r(G,A))
\rightarrow {\rm End}K_*(C^*_r(G,A))$$

The Baum-Connes conjecture for $G$ with coefficients in $A$ is therefore equivalent to the following statement:

 \noindent {\bf Conjecture}. The above map $R(G)\rightarrow {\rm End}K_*(C^*_r(G,A))$ sends the
Kasparov element $\gamma$ to the identity.

In this paper we prove the above conjecture  for the groups  $Sp(n,1)$ and $F_{4(-20)}$. The proof
also applies to the case of $SO(n,1)$ or $SU(n,1)$, but in that case one has in fact the much stronger result
$\gamma=1$ in $R(G)$: this was shown in [K2][JK], and also results from [HK], see for
example [J]. It also implies $K$-amenability (cf [JV]).

\bigskip
\noindent{\bf 1.2 Idea of our proof.}

\bigskip
The crucial point in our work is the role of {\it uniformly bounded representions}.
 Michael Cowling  has shown [C][ACD], already in 1980, that the non unitary principal series representions
$\pi _s$ on the suitable  Sobolev spaces 
associated to the hypoelliptic sublaplacian, are not only bounded, but uniformly bounded provided the complex
parameter lies in some vertical strip (section 3.1). 

On the other hand, there is a link between uniformly bounded representations and $K$-theory which we had
already observed in 1984 [J2]: the fact that the map $R(G)\rightarrow {\rm End}K_*(C^*_r(G,A))$ factors through a
group
$R_{\rm ub}(G)$ which is defined using uniformly bounded representations instead of unitary ones (section 3.2). 

The idea of our proof is therefore the following:
since we cannot (because of property $T$) make a homotopy from $\gamma$ to 1 with unitary representations, we
construct such a homotopy with uniformly bounded representations. The strip of uniformly bounded
representations of M. Cowling plays here the role of the complementary series for $S0(n,1)$ [K2] and $SU(n,1)$[JK].

Stricly speaking, there is another technical difficulty due to the fact that
at the  end of the homotopy  the representations are no more uniformly bounded: the uniform
bounds  blow up . This problem can be fixed using a trick which appears in
the Banach context of V. Lafforgue's thesis in 1998. The idea of adapting the argument to our situation had been
suggested to us in 1999 by N. Higson and V. Lafforgue. Technically we replace the uniformly bounded
representations  by representations of $\varepsilon$-exponential type, i.e. satisfying $\Vert\pi (g)\Vert\leq C
e^{\varepsilon l(g)}$, where $l(g)$ is the length of $g\in G$. We can perform a homotopy to 1 for any
$\varepsilon >0$ and a result of N. Higson and V. Lafforgue is that it is enough to show that our element of $R_{\rm
ub}(G)$ acts by the identity on the $K$-theory of any reduced crossed product $C^*_r(G,A)$.

\bigskip
\noindent{\bf 1.3. Uniformly bounded representations and $K$-theory.}
\bigskip

A representation $g\mapsto\pi (g)$ of a locally compact group $G$ in a Hilbert space is a strongly continuous
morphism of the group $G$ to the group of invertible elements of the $C^*$-algebra ${\cal L}(H)$ of bounded
operators on $H$.

\begin{definition}
  A $G$-Fredholm module is a triple $(H,\pi, T)$ where: $H$ is a ${\bf Z}/2$-graded Hilbert
space; $\pi$ a representation of $G$ in $H$ which is even  (i.e. commutes with the grading); $F$ a bounded operator
on $H$ which is odd (i.e. anticommutes with the grading), Fredholm (i.e. there exists a bounded operator $S$ on $H$
such that $TS-1$ and $ST-1$ are compact operators), and almost $G$-intertwining (i.e. $g\mapsto [T,\pi(g)]$ is a
normly continuous map from $G$ to the compact operators).
\end{definition}
\begin{definition}
 A homotopy of $G$-Fredholm modules is a triple $(E,\pi, T)$ where: $E$ is a ${\bf Z}/2$-graded
Hilbert module over the $C^*$-algebra $C[0,1]$ of continuous functions on $[0,1]$; $\pi$ a representation of $G$ in
$H$ which is even  (i.e. commutes with the grading);
$F$ a bounded operator on $E$ which is odd (i.e. anticommutes with the grading), Fredholm (i.e. there exists a
bounded operator
$S$ on $E$ such that $TS-1$ and $ST-1$ are compact operators), and almost $G$-intertwining (i.e. $g\mapsto
[T,\pi(g)]$ is a normly continuous map from $G$ to the compact operators).
\end{definition}

\begin{definition}
 A  bounded representation $\pi$ of $G$ is a representation of the group
$G$ by bounded operators on a Hilbert space which is strongly continuous. Such a
representation is uniformly bounded if moreover $$\sup_{g\in G}\Vert\pi (g)\Vert <\infty.$$
\end{definition}
Unitary representations are special cases of uniformly bounded representations.

\bigskip

The other crucial fact which makes the uniformly bounded representations relevant for our problem is the following.
Let $(G,A)$ be a $C^*$-dynamical system.  We will show below that the map $R(G)\rightarrow {\rm End}K_*(C^*_r(G,A))$
factors, as a map of additive groups,  through a group $R_{\rm u.b}(G)$ defined using uniformly bounded
representations instead of unitary representations.

We will need the following lemma:
\begin {lemma}
 1) Let $\pi$ be a uniformly bounded representation of $G$. Let $\lambda$ be
the left regular representation of $G$ on $L^2(G)$. There exists an operator $U$ on $H\otimes
L^2(G)$, which is bounded and has a bounded inverse, such that
$$\pi (g)\otimes\lambda (g)=U(1\otimes \lambda (g))U^{-1}$$
2) If moreover $\pi$ is a unitary representation, $U$ is a unitary operator.
\end{lemma}
The unitary case is classical and its generalisation to the uniformly bounded case is
straitforward.

Let $\lambda_A$ be the canonical map $C^*(G,A)\rightarrow C^*_r(G,A)$.

To a Hilbert space $H$ equipped with a uniformly bounded representation $\pi$,  associate the
Hilbert module $E=H\otimes C^*_r(G,A)$ and  the covariant representation of $(G,A)$ with
values in ${\cal L}_{C^*_r(G,A)}(E)$ defined by:

$$a\mapsto 1\otimes a,\  g\mapsto\pi(g)\otimes\lambda (g).$$ 

 \begin {lemma}
  The representation
 $\tilde\pi : L^1(G,A)\rightarrow {\cal L}_{C^*_r(G,A)}(E)$ extending the above covariant
representation factors through the reduced crossed product $C^*_r(G,A)$. 
\end{lemma}

\noindent {\it Proof.} The reduced crossed product  $C^*_r(G,A)$ is by
definition a sub-$C^*$-algebra of ${\cal L}_A(L^2(G)\otimes A)$ . Therefore the $C^*$-algebra
${\cal L}_{C^*_r(G,A)}(E)$ is a sub-$C^*$-algebra of ${\cal L}_A(H\otimes L^2(G)\otimes A)$.
It follows easily from lemma 1 that for any $a\in L^1(G,A)$, one has
$$\tilde \pi (a)=\tilde U(1\otimes\lambda_A(a))\tilde U^{-1}$$
where $$\lambda_A:L^1(G,A)\rightarrow C^*_r(G,A)\rightarrow{\cal L}_A(L^2(G)\otimes A)$$
The map $\tilde\pi$ therefore extends to a continuous (but not $*$) homomorphism
$$\tilde\pi: C^*_r(G,A)\rightarrow {\cal L}_{C^*_r(G,A)}(E).$$

 A $G$-Fredholm module (resp. a homotopy of $G$-Fredholm modules) is uniformly bounded if the representation $\pi$
is uniformly bounded.
 Let us consider the monoid of uniformly bounded $G$-Fredholm modules (for the direct sum),
and take the quotient by the equivalence relation defined by the uniformly bounded homotopies. We obtain an abelian
group denoted by $R_{\rm ub}(G)$. The same construction with "unitary" instead of "uniformly bounded" yields the
Kasparov group $R(G)=KK_G({\bf C},{\bf C})$ (we do not consider here the ring structure). Since any unitary map
is uniformly bounded, there is an obvious map $R(G)\rightarrow R_{\rm ub}(G)$.

 \begin{theorem}
 For any $G-C^*$-algebra $A$, the Kasparov map $R(G)\rightarrow {\rm End}K_*(C^*_r(G,A))$ factors
through the map $R(G)\rightarrow R_{\rm ub}(G)$.
\end{theorem}

\noindent{\it Proof.} Let us construct a map $R_{\rm ub}(G)\rightarrow {\rm End}K_*(C^*_r(G,A))$. To a $G$-Fredholm
module $(H,\pi ,T)$ we associate the triple $(H\otimes C^*_r(G,A), \tilde\pi, \tilde T)$ where 
$\tilde\pi: C^*_r(G,A)\rightarrow {\cal L}_{C^*_r(G,A)}(E)$ is the Banach algebra homomorphism defined above, and
$\tilde T=T\otimes 1\in {\cal L}_{C^*_r(G,A)}(H\otimes C^*_r(G,A))$.

 \begin{proposition} 
 $A$ is a Banach algebra, $B$ a $C^*$-algebra and $E$ a ${\bf
Z}/2$-graded $B$-Hilbert module, if $\pi$ is a banach algebra morphisme $A\rightarrow {\cal L}_B(E)$, and
$T\in{\cal L}_B(E)$ which is odd, Fredholm (i.e. invertible modulo the $B$-compact operators)  and commutes
modulo the $B$-compact operators
with the operators 
$\pi (a)$, $a\in A$. Then $(E,\pi, T)$ defines maps $K_0(A)\rightarrow K_0(B)$ and $K_1(A)\rightarrow K_1(B)$.
\end{proposition}

This is a classical statement. It is here unnecessary to use the powerful tool of Lafforgue's $KK^{ban}(A,B)$.

\bigskip
\noindent{\bf 1.4 Slow-growth representations.}

\bigskip
Let $\varepsilon >0$. We say that a representation $\pi$ of $G$ is of $\varepsilon$-exponential type if there is
a constant $C$ such that for any
$g\in G$, 
$$\Vert\pi (g)\Vert\leq C e^{\varepsilon l(g)}$$
where $l(g)=d(gx_0,x_0)$ is the length of $g$.
We define as above a $G$-Fredholm module (resp. a homotopy of $G$-fredholm modules of $\varepsilon$-exponential
type. Let $R_{\varepsilon}(G)$ be the abelian group of homotopy classes of such Fredholm modules.
There are obvious maps $R_{\varepsilon}(G)\rightarrow R_{\varepsilon '}(G)$ for $\varepsilon <\varepsilon '$, and we
consider the projective limit $\lim R_{\varepsilon}(G)$ when $\varepsilon\rightarrow 0$.
We would like to replace in the above result the group $R_{\rm ub}$ by $\lim R_{\varepsilon}(G)$. In fact there is a slightly weaker result due to Higson and Lafforgue (cf [L], th\'eor\`eme 2.3): 

 Theorem. The kernel of the map $$R(G)\rightarrow \lim R_{\varepsilon}(G)$$ is included in the kernel of the map $$R(G)\rightarrow {\rm
End}K_*(C^*_r(G,A))$$.

The idea of their proof is the following. For given $\varepsilon$ there is a Banach algebra $C_{\varepsilon}$ which is a completion of $C_c(G,A)$ such that the 
 the map
$$\tilde\pi : C_c(G,A)\rightarrow {\cal
L}_{C^*_r(G,A)}(E)$$
extending the  covariant representation of $(G,A)$ with
values in ${\cal L}_{C^*_r(G,A)}(E)$ defined by
$a\mapsto 1\otimes a,\  g\mapsto\pi(g)\otimes\lambda (g)$ extends to a map $C_{\varepsilon}\rightarrow{\cal L}_{C^*_r(G,A)}(E)$, which is compatible with the obvious maps $C_{\varepsilon}\rightarrow C_{\varepsilon'}$ for $\varepsilon'<\varepsilon$.

One has a commutative diagramme (cf [L] prop 2.5)
$$\matrix {R(G)&\rightarrow &{\rm
End}K_*(C^*_r(G,A))\cr\downarrow&&\downarrow\cr\lim R_{\varepsilon}(G)&\rightarrow&\lim{\rm Hom}(K_*(C_{\varepsilon}),K_*(C^*_r(G,A)))}$$

The theorem of Higson-Lafforgue follows immediatly from the lemma:
 Lemma. The group $K_*(C^*_r(G,A))$ is the union of the images of the maps $K_*(C_{\varepsilon})\rightarrow K_*(C^*_r(G,A))$.

The proof of their lemma is subtle. Using the idea of finite asymptotic dimension, they first show an estimate of the form (prop 2.6 in [L])
$$\Vert f\Vert_{C_{\varepsilon}}\leq k_{\varepsilon}e^{\varepsilon (ar+b)}\Vert f\Vert_{C^*_r(G,A)}$$
for $f\in C_c(G,A)$ with support in a ball of radius $r$ (for the length $l$).
Using the spectral radius formula in Banach algebras one immediately gets:
$$\rho_{C_{\varepsilon}}(f)\leq e^{\varepsilon ar}\rho_{C^*_r(G,A)}(f).$$
Hence the crucial fact that $\rho_{C^*_r(G,A)}(f)=\inf \rho_{C_{\varepsilon}}(f)$, which by standard holomorphic calculus implies the lemma.

\bigskip
\noindent{\bf 1.5 A characterisation of the element $\gamma$.}
\bigskip

We shall need a characterisation of the element  $\gamma$ of $R(G)$, which already appears  in [JK] . Consider $\bar
X=X\cup\partial X$  the natural compactification of $X$. 
We have the following fact:

\begin {proposition}
 The map $KK_G(C(\bar X),{\bf C})\rightarrow KK_K(C(\bar X),{\bf C})$ is an isomorphism.
In other words, $(1-\gamma )KK_G(C(\bar X),{\bf C})=0$.
\end{proposition}

The proof reproduces exactly the proof of [JK] , prop. 1.2.

 \begin{corollary}
  An element of R(G) which is in the image of $KK_G(C(\bar X),{\bf C})\rightarrow R(G)$ and maps
to 1 in $R(K)$ is equal to $\gamma$.
\end{corollary}

\bigskip \bigskip
\noindent {\bf 2. Construction of a family of Fredholm modules on the boundary of a symmetric space of rank 1.}

\bigskip
The main technical tools of the proof are introduced in thischapter . Let $X$ be a symmetric space of rank 1, and
$M=\partial X$ its boundary or sphere at infinity. Then $M$ is a generalized contact manifold, i.e. a smooth
manifold equipped with a subbundle $E$ of the tangent bundle. The group $G$ of isometries of $M$ leaves invariant
not only the generalized contact structure, but a conformal class of Carnot-Caratheodory metrics (i.e. metrics on
the subbundle
$E$ and on the quotient bundle $TM/E$) (section 2.1). 

To such a geometric data is associated a hypoelliptic sublaplacian and a pseudodifferiential calculus $\Psi
(M,E)$ [CGGP] modelled on the theory of representations a two step graded nilpotent groups (section 2.2).
This fact allows to construct the families of representations which play a central role in our proof. We start with
unitary principal representations, i.e. geometrically the spaces of $L^2$ sections of some equivariant
hermitian bundles on $M$. We deform them to non unitary principal series $\pi_s$ (depending on a complex parameter
$s\in{\bf C}$),
which by the virtue of the  pseudodifferential calculus $\Psi (M,E)$, become bounded on some Sobolev spaces 
associated to the hypoelliptic sublaplacian (section 2.3). On the other hand, there is a complex of differential
operators, introduced by M. Rumin [R],  which is more adapted than the de Rham complex to the generalized contact
and Carnot-Caratheodory structure and whose hypoelliptic regularity is expressed in terms of the pseudodifferential
calculus $\Psi (M,E)$ . A crucial fact that we prove in section 2.4 is the
$G$-invariance of the Rumin complex on the boundary of a rank one symmetric space.

All of the above geometric and analytic datas above are combined in the construction of a family (indexed by
$s\in{\bf C}$) of Fredholm modules (section 2.5).

\bigskip

\noindent{\bf 2.1 Field of graded Lie algebras associated to a Heisenberg manifold.}

\bigskip
We shall call Heisenberg manifold a smooth manifold $M$ equipped with a subbundle $E$ of the tangent bundle $TM$.  To such a subbundle is associated a bundle homomorphism
$${\bigwedge}^2E\rightarrow TM/E$$
which to $X\wedge Y$ associates $[X,Y]$ modulo $E$. This map is indeed well defined because the formula $[fX, gY]=fg[X,Y]+fX(g)Y-Y(g)fX\equiv fg[X,Y]$ modulo $E$ shows that the value of the vector field $[X,Y]$ modulo $E$ at a given point only depends on the value of $X$ and $Y$ at the same point.

 Note that the duality transpose of the above map is the map
$$E^{\perp}\rightarrow{\bigwedge}^2E^*$$
which to a form $\tau$ vanishing on the subbundle $E$  associates the restriction to $E$ of $d\tau$. Here again the formula $d(f\tau)=fd\tau+df\wedge\tau$ shows that the map is defined poinwise.

Let us consider $TM$ as a filtered bundle (by the increasing filtration $E\subset TM$) and let ${\rm gr}(TM)=E\oplus TM/E$ the associated graded bundle. Then there is a field of graded Lie algebra structures on the bundle ${\rm gr}(TM)$ defined (pointwise) by 
$$((X,Z),(Y,T))\mapsto (0,[X,Y]\ {\rm mod.}E)$$

The graded Lie algebra ${\rm gr}(TM)$ is equipped with the automorphism $\theta_t$ defined by $\theta_t (X)=tjX$ if $X\in  E$ and $\theta_t(Z)=t^2Z$ if $Z\in TM/E$

\noindent{\it Remark: } This construction can be generalized as follows. We consider a filtration of the tangent bundle $TM$:
$$F_1\subset F_2\subset ...\subset F_l=TM$$ 
subject to the following condition: if $X$ is a section of $F_i$ and $Y$ a section of $F_j$ then the Lie bracket $[X,Y]$ is a section of $F_{i+j}$.
Then we can equipp the associated graded bundle  ${\rm gr}(TM)=F_1\oplus F_2/F_1\oplus F_3/F_2... \oplus TM/F_{l-1}$
with a field of graded Lie algebra structure as follows. If $X$ is a section of $F_i$ and $Y$ a section of $F_j$, 
then the section $[X,Y]$ of $F_{i+j}$ modulo $F_{i+j-1}$ clearly only depends on $X$ modulo $F_{i-1}$ and $Y$ modulo $F_{j-1}$. 
Moreover the formula $[fX, gY]=fg[X,Y]+fX(g)Y-Y(g)fX$ again proves that the value of $[X,Y]\in F_{i+j}/F_{i+j-1}$ at a given point only depends on the values of
 $X\in F_i/F_{i-1}$ and  $Y\in F_j/F_{j-1}$ at the same point.
Note that the grading on the Lie algebra ${\rm gr}(TM)$ gives rise to the automorphism $\theta_t$ defined by $\theta_t (X)=t^jX$ if $X\in  F_j/F_{j-1}$.
 \bigskip\bigskip
 
 \noindent {\bf 2.2 Differential forms on a Heisenberg manifold.}
 
 \bigskip
The increasing filtration $E\subset TM$ defines a dual decreasing filtration $T^*M\supset E^{\perp}$, and induces filtrations of the exterior algebras $\bigwedge TM$ and 
$\bigwedge T^*M$. Namely:

$${\bigwedge}_1 TM\subset{\bigwedge}_2 TM\subset ...\subset\bigwedge TM$$
where $\bigwedge_wTM$ the space of multivectors of weight $\leq w$ is defined as the span of $X\wedge Y$ where $X\in\bigwedge^iE$, $Y\in\bigwedge^jTM$ and $i+2i\leq w$. Dually

$$\bigwedge T^*M\supset {\bigwedge}_1 T^*M\supset{\bigwedge}_2 T^*M\supset ...$$
where
$\bigwedge_w T^*M$, the space of forms of weight $\geq w$ is the span of
  $\alpha\wedge\beta$, $\alpha\in{\bigwedge}^iE^{\perp}$, $\beta\in{\bigwedge}^jT^*M$ where $i+2j\geq w$.
 Note that in the duality between $\bigwedge TM$ and $\bigwedge T^*M$, $\bigwedge_w T^*M$ is the orthogonal of $\bigwedge_{w-1} TM$.

 \begin{lemma}
  The associated graded algebra ${\rm gr} \bigwedge T^*M=\bigoplus_w\bigwedge_w T^*M/\bigwedge_{w+1} T^*M$ is canonically isomorphic to the exterior algebra $\bigwedge (E\oplus TM/E)^*$ of the dual of the Lie algebra ${\rm gr}(TM)=E\oplus TM/E$.
\end{lemma}

Let us consider the algebra $\Omega (M)=\Gamma (\bigwedge T^*M)$ of differential forms on $M$. Let $\Omega_w (M)$ be the space of sections of the bundle ${\bigwedge}_w T^*M$, 
defining a filtration of the algebra $\Omega (M)$. 
The associated graded algebra ${\rm gr}( \Omega (M))=\bigoplus\Omega_w(M)/\Omega_{w+1}(M)$ is the algebra of sections of the algebra bundle $\bigwedge (E\oplus TM/E)^*$.

 \begin{lemma}
  The map $d:Ê\Omega (M)\rightarrow\Omega (M)$ preserves the spaces $\Omega_w(M)$.  The induced map on the associated graded algebra is map $d_0$ given by the bundle map 
$$d_0:{ \bigwedge}^k(E\oplus TM/E)^*\rightarrow{\bigwedge}^{k+1}(E\oplus TM/E)^*$$
which is nothing but the coboundary map for the Lie algebra cohomology of the Lie algebra $E\oplus TM/E$.
\end{lemma}

\noindent{\it Remark: } This can be generalized as follows to the case of a filtration by the $F_j$'s. We define $\Lambda_w T^*M$ as the ideal of $\Lambda T^*M$ generated by the forms
  $\alpha_1\wedge\alpha_2\wedge ...\wedge\alpha_l$, with
$\alpha_j\in\Lambda^{i_j}F_{j-1}^{\perp}$ and the multiindex  $(i_1,...,i_l)$ satisfies  $i_1+2i_2+...+li_l=w$.
Both lemma 1 and lemma 2 can be generalized to the case of a filtration by $F_j$'s. The fact that $d$ preserves $\Omega_w(M)$ follows from the fact that $d(F_j^{\perp})\subset T^*M\wedge F_{j-1}^{\perp}$

We now consider the operator $d$ as a differential operator on $\Omega$. We equip the algebra ${\rm Diff}(M,\Lambda T^*M)$ of differential operators on the bundle $\Lambda T^*M$ with the following increasing filtration. Let us choose  a connection on $\Gamma (TM)$ such that  for any vector field $X$, $\nabla_X(\Gamma (E))\subset\Gamma (E)$. Such connections do exist.

We equipp he subalgebra of bundle operators (i.e. zero order operators) with the ${\bf Z}$-filtration such that:

(i) an element of ${\rm End}(\Lambda T^*M)$ is of weight $\leq w$ if it maps $\Lambda_u T^*M$ to $\Lambda_{u-w} T^*M$.

(ii)  if a vector field $X$ is in $\Gamma (E)$ (resp. $\Gamma (TM)$) then $\nabla_X$ has weight $\leq 1$ (resp.2). 

The filtration thus obtained is independent of the choice of the connection $\nabla$. 

 \begin{lemma}
 The operator $d$ in ${\rm Diff}(M,\Lambda T^*M)$ has 0-order in the above filtration.
\end{lemma}

The differential calculus given by ${\rm Diff}(M,\Lambda T^*M)$ with the above filtration extends to a pseudodifferential calculus $\Psi (M, E, \Lambda T^*M)$. Such a pseudodifferential calculus generalizes the pseudodifferential calculus associated to the tangent groupoid for the Heisenberg manifold [JvE], [vEY], see also [M].

\bigskip\bigskip

 \noindent{\bf 2.3 Quasi-conformal structure on a Heisenberg manifold.}
 
 \bigskip
 To any riemannian metric $\gamma$ on the manifold $M$ we  associate a bundle isomorphism $$\varphi_{\gamma}:{\rm gr}(TM)\rightarrow T^*M$$ as follows: we map ${\rm gr}(TM)$ isomorphically to $TM$ by identifying each quotient $F_j/F_{j-1}$ to the orthogonal (for the given metric) of $F_{j-1}$ in $F_j$, then the metric defines an isomorphism from $TM$ to its dual $T^*M$. Note that the map $\varphi_{\gamma}$ is filtration preserving, if ${\rm gr}(TM)$ is equipped with the decreasing filtration induced by its grading, and $T^*M$ by the decreasing filtration by the (duality) orthogonals to the increasing $F_j$'s.
 
  \begin{definition}
   Two metrics $\gamma_1$ and $\gamma_2$ are quasi-conformal if the composition $\varphi_{\gamma_1}^{-1}\circ\varphi_{\gamma_2}$ is a (filtered) Lie algebra isomorphism of ${\rm gr}(TM)$ which is 
 of the form 
 ${\rm Ad}(X)\circ\theta_t$ where $X\in{\rm gr}(TM)$ and $t\in{\bf R}_+^*$.
 \end{definition}
 
 The notion of quasi-conformality defines an equivalence relation among metrics on $M$. We shall call quasi-conformal structure on a Heisenberg manifold the choice of an equivalence class.  Given a quasi-conformal structure on a Heisenberg manifold there is a unique field of filtered Lie algebra structures on the cotangent bundle $T^*M$, such that the map $\varphi_{\gamma}:{\rm gr}(TM)\rightarrow T^*M$ is a filtered Lie algebra morphism (fiberwise) for any metric $g$ in the given class.
 
Let us assume that such a quasi-conformal class has been chosen, and let us consider the bundle $T^*M$ equipped by a field of filtered Lie algebra structures. Let $\delta :{ \bigwedge} ^kT^*M\rightarrow{\bigwedge} ^{k-1}T^*M$ be the bundle map which is on each fiber the boundary map for the Lie algebra homology of the Lie algebra $T^*M$. 

 \begin{lemma}
 The map $\delta$ preserves the filtration of  $\bigwedge T^*M$, and induces on the associated graded space $\bigwedge {\rm gr}(T^*M)$ the map $\delta_0=d_0^*$ , i.e. the adjoint of the coboundary map $d_0$ defined in lemma 2 for any metric $\gamma$ in the quasi-conformal class.
\end{lemma}

Let $\Omega =\Omega (M)$ be the graded algebra of differential forms on $M$. We consider on $\Omega$ the two operators $d$ and $\delta$, which  $d^2=0$ and $\delta^2=0$. They have respectively degree $1$ ans $-1$. 

Let us consider  the degree zero map $d\delta+\delta d$,  which commutes both with $d$ and $\delta$.

\begin{lemma}
  The map $d\delta+\delta d$ induces on the sections of the vector bundle ${\rm im}\delta$ a differential operator which is invertible among differential operators.
\end{lemma}

\noindent {\it Proof of the lemma.} Consider the filtration $\Omega_w (M)$ of $\Omega$. By lemmas 2 and 3 he maps $d$ and $\delta$ are both filtration preserving and induce on the associated graded space $\Gamma (\bigwedge {\rm gr}(T^*M))$ the map given by the bundle morphism $d_0d_0^*+d_0^*d_0$, which is clearly invertible on ${\rm im}d_0^*$. To deduce lemma 4, it is enough to observe that if a linear map on a filtered vector space induces an invertible map at the graded level, then it is invertible.
\bigskip\bigskip

\noindent{\bf 2.4 Complex associated to a quasi-conformal Heisenberg structure.}

\bigskip
Let $\cal E$ be the space of differential forms $\alpha$ on $M$ such that $\delta\alpha=0$ and $\delta d\alpha=0$. It is graded by $E^k={\cal E}\cap\Omega^k$ and is stable by $d$.  We thus have a subcomplex $({\cal E},d)$ of the complex $(\Omega, d)$.

Consider the map $\iota$ from $\cal E$ to the sections of the quotient bundle ${\rm ker}\delta/{\rm im}\delta$
obtained by composing the canonical injection
${\cal E}={\rm ker}\delta\cap{\rm ker}{\delta}d
 \rightarrow{\rm ker}\delta$ with the canonical surjection ${\rm ker}\delta\rightarrow{\rm ker}\delta/{\rm im}\delta$.
 Note that $({\rm ker}\delta/{\rm im}\delta)^k$ is the space of sections of the homology groups $H_k(\bigwedge T_x^*M, \delta )$
  of the Lie algebras $T_x^*M$ at each point $x\in M$.

  \begin{theorem}
  1) The  canonical injection ${\cal E}\rightarrow\Omega$ induces an isomorphism in cohomology. 
  2) The map  $\iota$ is an isomorphism from $\cal E$ to the sections of the quotient bundle ${\rm ker}\delta/{\rm im}\delta$.  . 
 \end{theorem}

{\it Proof of the theorem:}

Let us consider $q$ be the differential operator on $\Omega$ by $q\alpha= (d\delta+\delta d)^{-1}\delta\alpha$ where $(d\delta+\delta d)^{-1}$ is the inverse of $d\delta+\delta d$ on ${\rm im}\delta$. Note that the kernel of $q$ (resp. the image of $q$) is the space of sections of the vector bundle ${\rm ker}\delta$ (resp. ${\rm im}\delta$). 

 \begin {proposition}
One has $q^2=0$ and $qdq=q$.  The operator $\pi=dq+qd$ satisfies $\pi^2=\pi$, $\pi d=d\pi$ and $\pi q=q\pi=q$
\end{proposition}

Let us consider the decomposition $\Omega={\rm ker}\pi\oplus{\rm im}\pi$
It follows from lemma that
the cohomologies of the two subcomplexes $({\rm ker}\pi , d)$ and $({\rm im}\pi ,d)$ are:
 $$H^*({\rm im}\pi, d)=0, \    H^*({\rm ker}\pi, d)=H^*(\Omega, d).$$
 Indeed,  the map $dq+qd$ vanishes in  cohomology.
  
On the other hand  ${\rm ker}\pi={\rm ker}q\cap{\rm ker}qd$ and the map ${\rm ker}\pi\rightarrow{\rm ker}q/{\rm im}q$ is an isomorphism.
The first assertion follows immediately from the above lemma. The inverse of the last map is given by $1-qd$ which leaves ${\rm ker}q$ stable and vanishes on ${\rm im}q$

 \begin{proposition}
 The filtration $\Omega_w (M)$ of $\Omega$ restricts to a filtration of the space $\cal E$. It also induces a filtration of the bundle ${\rm ker}\delta/{\rm im}\delta$, and the map $\iota$ is filtration preserving. Moreover the filtration on ${\rm ker}\delta/{\rm im}\delta$ is associated to a canonical grading. So that the vector bundle ${\rm ker}\delta/{\rm im}\delta$ is given with a bigrading (degree and weight). 
\end{proposition}
\begin{corollary}
 The operator $D=\iota d{\iota}^{-1}$ defines a differential operator on the space of sections of the bundle ${\rm ker}\delta/{\rm im}\delta$ which satisfies $D^2=0$ and has degree one for the grading of ${\rm ker}\delta/{\rm im}\delta$. 
 \end{corollary}

The complex $\Gamma ({\rm ker}\delta/{\rm im}\delta),D)$ is called the BGG-Rumin complex. Its cohomology is the de Rham homology of the manifold $M$.

In particular the BGG-Rumin complex has finite dimensional cohomology. In the classical case (i.e. when the subbundle $E=TM$) this also follows from the ellipticity of the complex $(\Omega , d)$, namely the existence of a parametrix $Q_0$, a pseudodifferential operator of degree -1 such that $dQ_0+Q_0d-1$ is smoothing. Here we have to replace the classical pseudodifferential calculus by the calculus $\Psi (M, E, \Lambda T^*M)$ associated by Melin to a Lie filtration, here the filtration $E\subset TM$. 

 \begin{theorem}
 There exists an operator $Q_0\in\Psi (M, E, \Lambda T^*M)$ which has weight 0 (and degree -1 for the standard degree of forms) such that the operator $dQ_0+Q_0d-1$ is smoothing. 
\end{theorem}

To prove the theorem it is enough to consider the symbols, i.e. to look at the case of a nilpotent Lie group of Heisenberg type. For any $X$ of  degree 1  in the Lie algebra,  one has the Cartan formulas $di_X+i_Xd={\cal L}_X$ so that  $d$ and $i_X$ commute with ${\cal L}_X$. Let $(X_i)$ be an orthonormal basis of such vectores. Then $dA+Ad=\sum {\cal L}_{X_i}^2$ where $A=\sum i_{X_i}{\cal L}_{X_i}$, and $\sum {\cal L}_{X_i}^2$ can be shown to be invertible. Indeed ${\cal L}_X=X-{\rm ad}(X)^*$ and $\sum {\cal L}_{X_i}^2$ is equal to $\sum X_i^2$ plus a nilpotent operator. We then use the classical fact that $\sum X_i^2$ is invertible in the algebra of pseudodifferential operators on a Heisenberg type group.
\bigskip

We can restrict the operator $d$ to the subcomplex ${\cal E} $ of  $\Omega$ and consider the projection $\pi Q_0\pi$ as an operator on ${\cal E} $. By transport to $\Gamma ({\rm ker}\delta/{\rm im}\delta)$ we obtain a parametrix for the complex $D$. Let $\Psi (M, E, {\rm ker}\delta/{\rm im}\delta)$ be the corresponding pseudodifferential calculus.

 \begin{corollary}
  There exists an operator $Q_0\in\Psi (M, E, {\rm ker}\delta/{\rm im}\delta)$ which has weight 0 (and degree -1 for the standard degree of forms) such that the operator $DQ_0+Q_0D-1$ is smoothing.
\end{corollary}

\bigskip\bigskip
\noindent  {\bf 2.5 A family of representations.}
\bigskip

Let $(M,E)$ be a Heisenberg manifold equipped with a quasi-conformal structure manifold with $M$ compact, and $G$ be a Lie group acting on a 
$(M,E)$, preserving the structure. Let us fix  metric in the class and let
$\lambda_g$ be a cocycle as above. Let us consider the Hibert space $L^2(M)$ equipped with the unitary
representation $\pi$ defined by
$$\pi(g)f=\lambda_g^{\nu/2}g^{-1*}f$$

We shall use here the pseudodifferential calculus associated to the tangent groupoid for the Heisenberg manifold [JvE], [vEY], see also [M].

It follows from the above that for $g\in G$,
$\lambda_g^s\pi(g)L^s\pi(g)^{-1}-L^s$ is in $\Psi^{s-1}(M,E)$, and that 
$$L^{-s}\lambda_g^s\pi(g)L^s\pi(g)^{-1}-1$$
is in $\Psi^{-1}(M,E)$, and therefore compact.

 \begin{corollary}
  The operator $L^{-s}\lambda_g^s\pi(g)L^s$ extends to a bounded operator and 
$$L^{-s}\lambda_g^s\pi(g)L^s-\pi (g)$$ is compact. 
\end{corollary}

The map $g\mapsto L^{-s}\lambda_g^s\pi(g)L^s$ thus defines a representation of $G$ with values in bounded operators
on $L^2(M)$ with differs from the unitary representation $\pi (g)$ by compact operators.

\noindent {\it Remark.} The above construction can be generalized to the sections of a $G$-equivariant bundle $C$ equipped
with a hermitian metric $\Vert\ \Vert_C^2$ of weight $w\in {\bf R}$, i.e. such that $g^*\Vert\
\Vert_C^2=\lambda_g^{2w}$. One considers the Hilbert space $L^2(M,C)$ with the unitary representation defined by
$$\pi(g)f=\lambda_g^{\nu/2-w}g^{-1*}f$$
for any section $f$ of $C$. The operator $L$ is still defined by $L=(1+\Delta_E)^{1/2}$ , but now
$\Delta_E=\nabla_E^*\nabla_E$ where
$\nabla_E:C^\infty(M,C)\rightarrow C^*(M,C\otimes E^*)$ is the composition of a connection $\nabla :C^\infty
(M)\rightarrow C^\infty(M,C\otimes T^*M)$ with the restriction $T^*M\rightarrow E^*$ of 1-forms to the subbundle
$E$. The operator has the same analytical properties and the statement of the above corollary holds without any
change.
\bigskip

We shall need as a crucial point the result of M. Cowling. Let $M$ be as above the boundary of a
symmetric space of rank 1, with isometry group $G$, and $C$ be a
$G$-equivariant bundle  equipped with a hermitian metric $\Vert\ \Vert_C^2$ of weight $w\in {\bf R}$, i.e. such
that $g^*\Vert\
\Vert_C^2=\lambda_g^{2w}\Vert\
\Vert_C^2$. We equipp the Hilbert space $L^2(M,C)$ with the unitary representation defined by
$$\pi(g)f=\lambda_g^{\nu/2-w}g^{-1*}f$$
for any section $f$ of $C$. The operator $L$ is  defined by $L=(1+\Delta_E)^{1/2}$  as in section 2.3, remark.

 \begin{theorem}
 (M. Cowling) The representation $g\mapsto L^{-s}\lambda_g^s\pi(g)L^s$ extends to a bounded
operator for any $s\in{\bf C}$ such that $\vert{\rm Re}s\vert<\nu/2$. 
\end{theorem}

\bigskip\bigskip

\noindent {\bf 2.6 A family of Fredholm modules.}
\bigskip

   A $G$-Fredholm module is a triple $(H,\pi, T)$ where: $H$ is a ${\bf Z}/2$-graded Hilbert
space; $\pi$ a representation of $G$ in $H$ which is even  (i.e. commutes with the grading); $F$ a bounded operator
on $H$ which is odd (i.e. anticommutes with the grading), Fredholm (i.e. there exists a bounded operator $S$ on $H$
such that $TS-1$ and $ST-1$ are compact operators), and almost $G$-intertwining (i.e. $g\mapsto [T,\pi(g)]$ is a
normly continuous map from $G$ to the compact operators).

Note that in the above definition, the representation $g\mapsto\pi (g)$ does not have to be unitary.

 \begin{definition}
  A homotopy of $G$-Fredholm modules is a triple $(E,\pi, T)$ where: $E$ is a ${\bf Z}/2$-graded
Hilbert module over the $C^*$-algebra $C[0,1]$ of continuous functions on $[0,1]$; $\pi$ a representation of $G$ in
$H$ which is even  (i.e. commutes with the grading);
$F$ a bounded operator on $E$ which is odd (i.e. anticommutes with the grading), Fredholm (i.e. there exists a
bounded operator
$S$ on $E$ such that $TS-1$ and $ST-1$ are compact operators), and almost $G$-intertwining (i.e. $g\mapsto
[T,\pi(g)]$ is a normly continuous map from $G$ to the compact operators).
\end{definition}

From now on, $M=\partial X$ is the boundary of a symmetric space of rank 1, with isometry group $G$. We make a
choice of a point $x_0\in X$, which determines a Carnot-Caratheodory metric in the $G$-conformal class on the
generalized contact manifold $(M,E)$. We construct a family (indexed by $s\in{\bf C}$) of $G$-Fredholm modules
associated to the geometry of $M$

\begin{theorem}
 There exists a ${\bf Z}$-graded Hilbert space $H=\bigoplus_kH^k$, equipped with:
1) a family indexed
by
$s\in{\bf C}$ of representations $\rho_s$ of $G$ by bounded operators in $H$, of degree zero for the ${\bf
Z}$-grading, such that $\rho_0$ is unitary and $\rho_s(g)-\rho_0(g)$ is compact for any $s\in{\bf C}$ and any $g\in
G$.
2) a bounded operator $F$ in $H$, of degree one for the ${\bf Z}$-grading, with $F^2=0$, and such
that$[F,\rho_s(g)]$ is compact for any $s\in{\bf C}$ and $g\in G$, and $[F,\rho_1(g)]=0$.
3) a bounded operator $Q$ in $H$, of degree $-1$ for the ${\bf Z}$-grading, such that the operators $Q^2$, $FQ+QF-1$
and $[Q,\rho_s(g)]$ are compact for any $s\in{\bf C}$ and $g\in G$.
\end{theorem}
In particular, one has

  \begin{corollary}
   For any $s\in {\bf C}$, the triple $(H,\rho_s , F+Q)$ is a $G$-Fredholm module. The Fredholm
modules for different values of $s$ are homotopic.
\end{corollary}

Let us construct the objects appearing in the theorem.

i) Let $H$ be the Hilbert space of $L^2$ sections of the bundle $C$ on $M=\partial X$. One has a bigrading 
$H=\bigoplus H^{i,j}$ where the sum is over $0\leq i\leq p$ and $0\leq j\leq q$, and $H^{i,j}$ is the space of $L^2$
sections of $C^{i,j}$. We shall consider $H$ as  ${\bf Z}$- graded as follows: $H=\bigoplus_kH^k$ where
$$H^k=\bigoplus_{i+j=k}H^{i,j}.$$ 
The ${\bf Z}$-grading induces as usual a ${\bf Z}/2$-grading, which grading operator
$\gamma=(-1)^{i+j}$ on
$H^{i,j}$.

ii) For any $s\in{\bf C}$ , let $\pi_s$ be the representation of $G$ by bounded operators on each Hilbert space
$H^{i,j}$ defined by 
$$\pi_s(g)\alpha =\lambda_g^{(\nu/2-w)(1-s)}g^{-1*}\alpha$$
for any section $\alpha$ of $C^{i,j}$, where $w=i+2j$.
Note that $\pi_s(g)=\lambda_g^{-(\nu/2-w)s}\pi_0(g)$ where $\pi_0$ is a unitary representation.

We also define $$\rho_s(g)=L^{(\nu/2-w)s}\pi_s(g)L^{-(\nu/2-w)s}$$
where again $w=i+2j$.
By section 2.3, $\rho_s$ is a representation of $G$ by bounded operators on $H^{i,j}$, and $\rho_s(g)-\pi_0(g)$ is
compact.

We shall also denote $\pi_s$ and $\rho_s$ the representations of $G$ on the Hilbert space $H$ obtained by direct
sum over all the $(i,j)$'s. One has 
$$\rho_s(g)=L^{(\nu/2-W)s}\pi_s(g)L^{-(\nu/2-W)s}$$
where $L^{(\nu/2-W)s}$ is the non homogeneous pseudodifferential operator on sections of $C$ defined by
$L^{(\nu/2-w)s}$ on sections of $C^{i,j}$ with $i+2j=w$. The operators $\rho_S(g)$ are of degree zero for
the ${\bf Z}$-grading, and thus  even operators for the
${\bf Z}/2$-grading.

Note that the parametrisation has been chosen such that for $s=1$, the representation $\pi_1$ is the natural
representation of $G$ acting functorially on the sections of the $G$-equivariant bundle $C$.

iii) Let $F$ be the Rumin operator made bounded as follows:
$$F=L^{\nu/2-W}DL^{-(\nu/2-W)}.$$
By section 2.4, $F$ is bounded since it belongs to $\Psi^0(M,E)$. It is  of degree $1$ for the
${\bf Z}$-grading, and thus an odd operator for the
${\bf Z}/2$-grading. The $G$-invariance of $D$ means that
$$[D,\pi_1(g)]=0,$$
which implies 
$$[F,\rho_1(g)]=0.$$
But the fact that $\rho_s(g)-\pi_0(g)$ is compact for any $s$ , together with the boundedness of $F$ implies that
$$[F,\rho_s(g)]$$
is compact for any $g\in G$ and $s\in{\bf C}$.

 \begin{lemma} 
 There exists a $Q\in\Psi^0(M,E)$ such that $FQ+QF-1$ is a smoothing operator. In particular $Q$ is
bounded and $FQ+QF-1$ is compact.
Moreover, one has gor $g\in G$, $\rho_1(g)Q\rho_1(g)^{-1}-Q\in \Psi^{-1}(M,E) $, in particular $[Q,\rho_1(g)]$ is
compact  and so is $[Q,\rho_s(g)]$ for any $s$.
\end{lemma}

The lemma follows from the corollary in 2.4  by considering 
$$Q=L^{\nu/2-W}Q_0L^{-(\nu/2-W)}$$.
 
 \noindent{\it Remark:} one can always assume that $Q^2$ is smoothing. Namely il is easy to check that $FQ$ and $QF$ are idempotent modulo smoothing, and that $FQ^2F$ is smothing. Now, replacing $Q$ by $\tilde Q=QFQ$ we note that $\tilde Q-Q$ is smoothing, that $\tilde Q^2$ is smothing and satisfies $F\tilde Q+\tilde Q F=1$ modulo smoothing.

 When $s=0$, the representation $\rho_0$ is unitary, so that our Fredholm module is a representative of a class in
Kasparov's ring $R(G)=KK_G({\bf C},{\bf C})$. When
$s=1$, the representation $\rho_1$ exactly commutes with $F$, a fact which implies the Fredholm module is homotopic
to a trivial module since $F+Q$ can be replaced by $t^{-1}F+tQ$ whose commutator with $\rho_1(g)$ vanishes when
$t\rightarrow 0$.

Note also that the index of the Fredholm operator $F+Q$ is equal to the Euler caracteristic of the
complex $D$, i.e. the Euler caracteristic of the manifold $M$ : it is $2$ or $0$ depending whether the sphere $M$
is even ( the $g=SO(2n+1,1)$-case) or odd dimensional (all other cases).  We shall have to  modify the familily of
Fredholm modules  to get an index one operator. {\it This modification will be explained in section 3É}

\bigskip\bigskip
\noindent {\bf 2.7. Uniform boundedness.} 
\bigskip
 
The following crucial fact follows from M. Cowling's theorem stated at the end of 2.5.

 \begin{proposition}
 The representation $\rho_s$ on $H^{i,j}$ is uniformly bounded for any $s\in{\bf C}$ such that
$$\vert{\rm Re}s\vert<{\nu/2\over\vert\nu/2-w\vert}$$
with $w=i+2j$.
\end{proposition}

Note that the values of $w$ are between $0$ and $\nu=p+2q$, so that $\nu/2\leq \vert\nu/2-w\vert$ with equality iff
$w=0$ or $\nu$, i.e. $(i,j)=(0,0)$ or $(p,q)$. In particular, 

 \begin{corollary}
  The representations $\rho_s$ are uniformly bounded if $s\in [0,1[$; the representation
$\rho_1$ is uniformly bounded on $H^{i,j}$ for $(i,j)\neq(0,0)$ and $\neq (p,q)$.
\end{corollary}

 A $G$-Fredholm module (resp. a homotopy of $G$-Fredholm modules) is uniformly bounded if the representation $\pi$
is uniformly bounded. The $G$-Fredholm module $(H,\rho_s , F+Q)$ is uniformly bounded for any $s\in [0,1[$

\bigskip\bigskip 
\noindent{\bf 3. Sketch of proof of the Baum-Connes conjecture for $Sp(n,1)$ and $F_{4(-20)}$.}
\bigskip

We shall first modify the  Fredholm module $(H,\rho_s , F+Q)$ to make it of index 1. We will then obtain an element
$\delta$ of $R_{\rm ub}(G)$ and show that it is a representative of the Kasparov element $\gamma\in R(G)$. We would
like to show that it is  equal to 1 in the group
$R_{\rm ub}(G)$, but in fact we show a weaker statement, namely that it is equal to 1 in $R_{\varepsilon}(G)$ for
any small $\varepsilon$ . That is enough to assure that it acts by the identity on the
$K$-theory of reduced crossed products, therefore proving the Baum-Connes conjecture (with coefficients).

\bigskip\bigskip
\noindent {\bf 3.1. The trunkated Fredholm module.}
\bigskip

To obtain a $G$-Fredholm module of index one, we modify the triple $(H,\rho_s,F+Q)$ as follows. The idea is to
trunkate the complex in the middle.

1) Let $H'=\bigoplus_kH'^k$ be the ${\bf Z}$-graded Hilbert space defined as follows. 

For $k\leq {m-1\over 2}$, we let $H'^k=H^k$.

For $k={m+1\over 2}$, we define $H'^{m+1\over 2}=F(H^{m-1\over 2})$

For $k>{m+1\over 2}$, we let $H'^k=0$. 

2) The representation $\rho'_s$ is defined as follows.
For $k\leq {m-1\over 2}$, we let $\rho'_s=\rho_s$ on $H'^k=H^k$

For $k={m+1\over 2}$, we define $\rho'_s=\rho_1$ (independent on $s$) which stabilizes $H'^{m+1\over
2}=F(H^{m-1\over 2})$ because $[F,\rho_1(g)]=0$.
Note that we still have $\rho'_s(g)-\rho_0(g)$ compact. The representation $\rho'_s$ is uniformly bounded for
$s\in [0,1[$. 

3) The operator $F'$ is defined as follows: $F':H'^k\rightarrow H'^{k+1}$ is the restriction of $F$ to $H^k$ for
any $k\leq {m-1\over 2}$. Its parametrix $Q'$ is defined similarly: $Q':H'^k\rightarrow H^{k-1}$ is the restriction
of $Q$ to $H'^k$ for all $k\leq {m+1\over 2}$.

The triple $(H',\rho'_s,F'+Q')$ is a Fredholm module, with index one for any $s$. It is uniformly bounded for
$s\in [0,1[$. Let $\delta$ be its class in   $R_{\rm u.b}(G)$, which is independent of $s$.
\bigskip

\noindent{\bf 3.2 Representing the $\gamma$ element.}

\bigskip
We show that the trunkated Fredholm module defined in the previous section is a representative of Kasparov's
$\gamma$ element of $R(G)$. More precisely:

 \begin{theorem} 
 The class $\delta$ is the image of the Kasparov element $\gamma$ under the map
$R(G)\rightarrow R_{\rm u.b.}(G)$.
\end{theorem}

\noindent {\it Proof.} We take $s=0$. The representation  on all $H'^k$  for $k\neq {m+1\over 2}$ is $\rho_0$
which is unitary  , but it is $\rho_1$ on $H'^{m+1\over 2}$ so that it is only uniformly bounded. The key lemma is
that $\rho_1$ on
$H'^{m+1\over 2}$ is in fact equivalent to a unitary representation.

The crucial lemma is the following:

 \begin{lemma} There exists a  bounded operator $U$ from $H'^{m+1\over 2}$ to the Hilbert space 
${\cal H}^{m+1\over 2}$ of
$L^2$-harmonic  forms (of degree ${m+1\over 2}$) on the symmetric space $X$, which is invertible with bounded
inverse, and such that
$U\rho_1(g)=\sigma(g)U$ for any $g\in G$, where $\sigma$ is the natural unitary representation of $G$ on 
${\cal H}^{m+1\over 2}$ defined by $\sigma (g)\omega =g^{-1*}\omega$.
\end{lemma}

 \begin{corollary}
  The element $\delta$ is in the image of the map  $R(G)\rightarrow
R_{\rm ub}(G)$.
\end{corollary}

\noindent {\it Proof of the corollary.} Let us show that the above Fredholm module $(H', \rho'_0, F+Q)$ is conjugate
to a unitary Fredholm module. We define $$V:H'=\bigoplus_{k=0}^{m-1\over 2}H^k\oplus H'^{m+1\over 2}
\rightarrow \bar H=\bigoplus_{k=0}^{m-1\over 2}H^k\oplus{\cal H}^{m+1\over 2}$$
by $V=1\oplus U$. Then the conjugate Fredholm module $(\bar H, V\rho'_0(.)V^{-1}, V(F+Q)V^{-1})$ , defining the
same class in $R_{\rm ub}(G)$, is unitary since $V\rho'_0(.)V^{-1}=\rho_0\oplus U\rho_1(.)U^{-1}=\rho_0\oplus\sigma$
is a unitary representation.
\bigskip

To show that $\delta$ is actually the image of $\gamma$, one shows as in [JK] that it is in the image of 
$KK_G(C(\bar X),{\bf C})$ where $\bar X=X\cup\partial X$ is the compactification of $X$. Indeed, the spaces $H^k$
of $L^2$-sections of bundles on $\partial X$ carry  representations of the commutative $C^*$-algebra $C(\partial
X)$ and hence of $C(\bar X)\rightarrow C(\partial X)$ by restriction. On the other hand, $C(\bar X)$ is represented
in the Hilbert space $L^2\Omega (X)=L^2(X,\bigwedge T^*X)$ of $L^2$-forms on $X$. Let us replace in $\delta$ the spaces 
$\bar H^k=H^k$ ($0\leq k\leq {m-1\over 2}$) and $\bar H^{m+1\over 2}={\cal H}^{m+1\over 2}$ by 

$$\tilde H^k=H^k\oplus L^2\Omega^k(X)$$ 

for $0\leq k\leq {m-1\over 2}$
and $$\tilde H^k=L^2\Omega^k(X)$$ for ${m+1\over 2}\leq k\leq m+1$

Note that $$\tilde H^k=H'^k\oplus{{\cal H}^k}^{\perp}$$ for all $k$, where of course the space of harmonic $L^2$-forms ${{\cal H}^k}$ is non zero only in degree $k={m+1\over 2}$.
Let us extend the operator $F'$ by the direct sum with the phase of the de Rham operator $d$ on $L^2\Omega^k(X)$. Then one shows as in [JK] that we still get a representative of $\delta$. The key observation is the following lemma:

 \begin{lemma} The above map $UF:H'^{m-1\over 2}\rightarrow L^2(X,\bigwedge^{m+1/2}T^*X)$ commutes  modul o
compact operators with the action of $C(\bar X)$.
\end{lemma}

To summarize we obtain:

\begin{corollary} The element $\delta$ is in the image of the composition of maps 
$$KK_G(C(\bar X),{\bf C})\rightarrow R(G)\rightarrow R_{\rm ub}(G).$$
\end{corollary}

The proposition in 1.5 (cf. [JK] , prop. 1.2.) then implies  theorem 1.

\bigskip\bigskip
\noindent {\bf 3.3 Homotopy to 1.}
\bigskip

 \begin{theorem}
  For any $\varepsilon >0$, the class $\delta$ maps to $1$ under the map $$R_{\rm
u.b}(G)\rightarrow R_{\varepsilon}(G).$$
\end{theorem}
\begin{corollary}
 For any $G-C^*$-algebra $A$, the class $\delta$ maps to the identity under the map 
$$R_{\rm ub}(G)\rightarrow {\rm End}K_*(C^*_r(G,A)).$$
\end{corollary}

 The above theorem relies on the following

 \begin{lemma}
 Let $\varepsilon \in{\bf R}$, $s\in {\bf C}$ and $g\in G$. The operator
$$\rho_s^{\varepsilon}(g)=L^{-\varepsilon}\rho_s(g)L^{\varepsilon}$$ on $H^{i,j}$, $i+2j=w<\nu$
is bounded for any $g\in G$. Assume furthermore that $\varepsilon >0$ if $w=0$ and $\varepsilon <0$ if $w=\nu$. The
map
$g\mapsto
\rho_s^{\varepsilon}(g)$ defines a representation of
$G$ satisfying  
$\Vert\rho_s^{\varepsilon}(g)\Vert\leq Ce^{\varepsilon al(g)}$, with  constants $C$ and $a$ independant of $g$
and $s$.
\end{lemma}

\bigskip
\noindent {\it Proof of the lemma.} One has
$\pi_{s-\varepsilon}(g)=\lambda_g^{\varepsilon (\nu/2-w)}\pi_s(g)$ , therefore
$$\rho_s^{\varepsilon}(g)=L^{(\nu/2-w)s-\varepsilon}\lambda_g^{-\varepsilon} L^{-(\nu/2-w)s+\varepsilon
}\rho_{s-{\varepsilon\over \nu/2-w}}(g).$$
where the operator $L^{\nu/2-w-\varepsilon}\lambda_g^{-\varepsilon} L^{-(\nu/2-w-\varepsilon
)}$ is in $\Psi^0(M,E)$, and therefore bounded.

Note that the representations  $\rho_{s-{\varepsilon\over \nu/2-w}}(g)$ are uniformly bounded with a bound
independant of $s$ because $s-{\varepsilon\over \nu/2-w}$ belongs to a compact interval strictly contained in
$[0,{\nu/2\over\vert\nu/2-w\vert}[$. Moreover, the map $g\mapsto L^{\nu/2-w-\varepsilon}\lambda_g
L^{-(\nu/2-w-\varepsilon )}$ is a 1-cocycle with respect to the representation $\rho_{s-{\varepsilon\over
\nu/2-w}}(g)$ and satisfies an inequatility 
$$\Vert L^{\nu/2-w-\varepsilon}\lambda_g
L^{-(\nu/2-w-\varepsilon )}\Vert\leq C e^{al(g)}$$
with positive constants $C$ and $a$. 

\bigskip
\noindent{\it Proof of  theorem 2.}

 Let us consider 
$F_{\varepsilon }=L^{\varepsilon}FL^{-\varepsilon}$ and $Q_{\varepsilon }=L^{\varepsilon}QL^{-\varepsilon}$.
Then $[F^{\varepsilon }, \rho_s^{\varepsilon}(g)]$ is compact for any $g\in G$ . The triple
$(H,\rho_s^{\varepsilon}, F_{\varepsilon}+Q_{\varepsilon})$ is a $G$-Fredholm module of $\varepsilon$-exponential
type for any $s\in[0,1]$, which represents the same homotopy class
in
$R_{\varepsilon}(G)$. Now for $s=1$, 
 $[F^{\varepsilon }, \rho_1^{\varepsilon}(g)]=0$ so that the Fredholm module is trivial, i.e. equal to its index
which is 1. 

\bigskip\bigskip\bigskip

\centerline{\bf Bibliography.}

\bigskip\bigskip

[ACD] F. Astengo, M. Cowling, B. Di Blasio, preprint 2001.

[BCH] P. Baum, A. Connes, N. Higson, {\it Classifying space for proper actions and $K$-theory
of group $C^*$-algebras}, in: $C^*$-algebras: 1943-1993, Contemp. Math. 167, AMS, 1994,
240-291.

[C] M. Cowling, {\it Unitary and uniformly bounded representations of some simple Lie
groups},in :
Harmonic analysis and group representations, 49--128, Liguori, Naples, 1982. 

[CGGP] M. Christ, D. Geller, P. Glowacki, L. Polin, {\it Pseudodifferential operators on
groups with dilatations}, Duke Math. J. 68 (1992), 31-65.

[HK] N. Higson, G. Kasparov, {\it $E$-theory and $KK$-theory for groups which act properly and
isometrically on Hilbert space}, Invent. Math. 144 (2001), no. 1, 23--74. 

[HL] N. Higson, V. Lafforgue, preprint 1999.

[J1]  P. Julg, {\it Complexe de Rumin, suite spectrale de Forman et cohomologie $L^2$ des
espaces sym\'etriques de rang $1$}, C. R. Acad. Sci. Paris S\'er. I Math. 320 (1995),
no. 4, 451--456. 

[J2] P. Julg, {\it Remarks on the Baum-Connes conjecture and Kazhdan's property $T$}, Fields
Instit. Commun., volume 13, Amer. Math. Soc., 1997, 145-153.

[J3] P. Julg, {\it Travaux de Higson et Kasparov sur la conjecture de Baum-Connes},
s\'eminaire Bourbaki, expos\'e 841, mars 1998, Ast\'erisque 252, 151-183.

[J4] P. Julg, {\it La conjecture de Baum-Connes \`a coefficients pour le groupe $Sp(n,1)$},  C. R. Acad.
Sci. Paris , s\'er.I, 334 (2002), 533-538.

[JK] P. Julg, G. Kasparov, {\it Operator $K$-theory for the group $SU(n,1)$}, J. Reine
Angew. Math. 463 (1995), 99-152.

[JV] P. Julg, A. Valette, {\it $K$-theoretic amenability for 
${\rm SL}_2({\bf Q}_p)$, and the action on the associated tree}, J. Funct. Anal. 58 (1984), no. 2,
194--215.

[K1] G.G. Kasparov, {\it $K$-theory, group $C^*$-algebras and higher signatures (Conspectus}, preprint
Chernogolovka 1981, reprinted in {\it Novikov Conjectures, Index Theorems and Rigidity}, edited by S. Ferry, A.
ranicki and J. Rosenberg, London Mathematical Society Lecture Note Series 226, Cambridge 1995; vol. 1, 101-146.

 [K2] G. Kasparov, {\it Lorentz groupes: $K$-theory of unitary representations
and crossed products}, Soviet. Math. Dokl. 29 (1984), n. 2, 256-260.

[K3] G. Kasparov,{\it  Equivariant $KK$-theory and the Novikov conjecture}, Invent. Math. 91
(1988), n. 1,147-201.

[L] V. Lafforgue, {\it Une d\'emonstration de la conjecture de Baum-Connes pour les groupes
r\'eductifs sur un corps $p$-adique et pour certains groupes discrets poss\'edant la
propri\'et\'e $(T)$}, C. R. Acad. Sci. Paris S\'er. I
Math. 327 (1998), no. 5, 439--444. 

[R] M. Rumin,{\it Differential geometry on $C-C$ spaces and application to the Novikov-Shubin
numbers of nilpotent Lie groups}, C. R. Acad. Sci. Paris S\'er. I Math. 329 (1999), no. 11,
985--990. 

[W] A. Wassermann, {\it Une d\'emonstration de la conjecture de Connes-Kasparov pour les
groupes de Lie lin\'eaires connexes r\'eductifs}, C. R. Acad. Sci. Paris S\'er. I Math. 304
(1987), no. 18, 559--562.

\bigskip\bigskip

Pierre Julg

Universit\'e d'Orl\'eans, MAPMO et FDP

BP 6759

F-45067 Orl\'eans Cedex 2

pierre.julg@univ-orleans.fr

\end{document}